\documentclass[a4paper]{article}
\usepackage{latexsym, amssymb, amsthm}
\begin{document}
\theoremstyle{plain}
\newtheorem{theorem}{Theorem}[section]

\newtheorem{lemma}[theorem]{Lemma}
\newtheorem*{Cor}{Corollary}
\newtheorem{corollary}[theorem]{Corollary}
\mathsurround 2pt

\title{On products of groups with abelian subgroups of small index}
\author{Bernhard Amberg and Yaroslav Sysak\thanks
{The second author likes to thank the Deutsche Forschungsgemeinschaft 
for financial support and the Institute of Mathematics of the University of Mainz 
for its excellent hospitality during the preparation of this paper.}}
\renewcommand{\thefootnote}{\fnsymbol{footnote}}
\setlength{\footnotesep}{.1in} 
\date{}

\bibliographystyle{plain}

\maketitle
\begin{abstract}
It is proved that every group of the form $G=AB$ with two subgroups $A$ and $B$ each of which is either abelian or has a quasicyclic subgroup of index $2$ is soluble of derived length at most $3$. In particular, if $A$ is abelian and $B$ is a locally quaternion group, this gives a positive answer to Question 18.95 of "Kourovka notebook" posed by A.I.Sozutov. \end{abstract}

\section{Introduction}

Let the group $G=AB$ be the product of two subgroups $A$ and $B$, i.e. $G=\{ab\mid a\in A, b\in B\}$. It was proved by N. It\^ {o} that the group $G$ is metabelian if the subgroups $A$ and $B$ are abelian (see \cite[Theorem 2.1.1]{AFG}). 

In connection with  It\^ {o}'s theorem a natural question is whether every group $G=AB$ with abelian-by-finite subgroups $A$ and $B$ is metabelian-by-finite \cite[Question 3]{AFG} or at least soluble-by-finite. However, this seemingly simple question is very difficult to attack and only partial results in this direction are known. A positive answer was given for linear groups $G$ by the second author in \cite{Sy1} (see also \cite{S_09}) and for residually finite groups $G$  by J. Wilson \cite[Theorem 2.3.4]{AFG}. Furthermore,  N.S. Chernikov proved that every group $G=AB$ with central-by-finite subgroups $A$ and $B$ is soluble-by-finite (see \cite[Theorem 2.2.5]{AFG}). 

It is natural to consider first groups $G=AB$ where the two factors $A$ and $B$ have abelian subgroups with small index. There are a few known results in the case when both factors $A$ and $B$ have an abelian subgroup of index at most $2$. It was shown in \cite{AS_08} that $G$ is soluble and metacyclic-by-finite if $A$ and $B$ have cyclic subgroups of index at most $2$, and it is proved in \cite{AFK} that $G$ is soluble if $A$ and $B$ are periodic locally dihedral subgroups. A more general result that  $G=AB $ is soluble if each of the factors $A$ and $B$ is either abelian or generalized dihedral was obtained in \cite{AS_13} by another approach. Here a group is called generalized dihedral if it contains an abelian subgroup of index $2$ and an involution which inverts the elements of this subgroup. Clearly dihedral groups and  locally dihedral groups, i.e. groups with a local system of dihedral subgroups, are generalized dihedral.   

We recall that a group is called quasicyclic (or a Pr\" ufer group) if it is an infinite locally cyclic $p$-group for some prime $p$. It is well-known that quasicyclic subgroups of abelian groups are their direct factors. Furthermore, it seems to be known and will be shown below that every non-abelian group having a quasicyclic subgroup of index $2$ is either an infinite locally dihedral or a locally quaternion group. It should be noted that for each prime $p$, up to isomorphism, there exists a unique locally dihedral group whose quasicyclic subgroup is a $p$-group, and there is only one locally quaternion group. These and other details about such groups can be found in \cite{KW}, p. 45 - 50. 

\begin{theorem}\label{Main} Let the group $G=AB$ be the product of two subgroups $A$ and $B$ each of which is either abelian or has a quasicyclic subgroup of index $2$.  Then $G$ is soluble with derived length at most $3$. Moreover, if  the subgroup $B$ is non-abelian and $X$ is its quasicyclic subgroup, then $AX=XA$ is a metabelian subgroup of index $2$ in $G$. \end{theorem}

As a direct consequence of this theorem,  we have an affirmative answer to Question 18.95 of the "Kourovka notebook" posed by A.I.Sozutov.

\begin{corollary}\label{1} If a group $G=AB$ is the product of an abelian subgroup $A$ and a locally quaternion subgroup $B$, then $G$ is soluble. \end{corollary}

It is also easy to see that if each of the factors $A$ and $B$ in Theorem \ref{Main} has a quasicyclic subgroup of index $2$, then their quasicyclic subgroups are permutable.  As a result of this the following holds.

\begin{corollary}\label{2} Let the group $G=A_1A_2\cdots A_n$ be the product of pairwise permutable subgroups $A_1$, ..., $A_n$  each of which contains a quasicyclic subgroup of index $2$.  Then the derived subgroup $G'$ is a direct product of the quasicyclic subgroups and the factor group $G/G'$ is elementary abelian of order $2^m$ for some positive integer $m\le n$. \end{corollary}

The notation is standard. If $H$ is a subgroup of a group $G$ and $g\in G$, then the normal closure of $H$ in $G$ is the normal subgroup of $G$ generated by all conjugates of $H$ in $G$, and  $g^G$ is  the conjugacy class of $G$ containing $g$, respectively.

\section{Preliminary lemmas}

Our first lemma lists some simple facts concerning groups with quasicyclic subgroups of index $2$ which will be used without further explanation. 

\begin{lemma}\label{qua} Let $G$ be a non-abelian group containing a quasicyclic $p$-subgroup $X$ of index $2$ and $y\in G\setminus X$. Then $y^2\in X$ and the following statements hold:\begin{itemize}
\item[1)] every subgroup of $X$ is characteristic in $G$; 
\item[2)] the group $G$ is either locally dihedral or locally quaternion;
\item[3)] the derived subgroup $G'$ coincides with $X$;
\item[4)] every proper normal subgroup of $G$ is contained in $X$;
\item[5)] if $G$ is locally quaternion, then $p=2$, $y^4=1$, $x^y=x^{-1}$ for all $x\in X$,  
the center $Z(G)$ coincides with $\langle y^2\rangle$ and is contained in every non-trivial subgroup of $G$, 
the coset $yX$ coincides with the conjugacy class $y^G=y^X$;
\item[6)] if $G$ is locally dihedral, then $y^2=1$, $x^y=x^{-1}$ for all $x\in X$,  $Z(G)=1$ and the coset $yX$ coincides with the conjugacy class $y^G=y^X$ for $p>2$ and $Z(G)$ is the subgroup 
of order $2$ in $X$ for $p=2$;   
\item[7)] the factor group $G/Z(G)$ is locally dihedral.
\end{itemize}
\end{lemma} 

\proof In fact, only statement 2) needs an explanation.  Clearly $G=X\langle y\rangle$ for some $y\in G$ with $y^2\in X$ and each cyclic subgroup $\langle x\rangle$ of $X$ is normal in $G$.  Therefore for $p > 2$ we have $y^2=1$ and either $x^y=x$ or $x^y=x^{-1}$. Since $X$ contains a unique cyclic subgroup of order $p^n$ for each $n\ge 1$, the equality $x^y=x$ for some $x\ne 1$ holds for all $x\in X$, contrary to the hypothesis that $G$ is non-abelian. Therefore $x^y=x^{-1}$ for all $x\in X$ and hence the group $G$ is locally dihedral. In the case $p=2$  each subgroup $\langle x\rangle$ of $X$ properly containing the subgroup $\langle y^2\rangle$ has index $2$ in the subgroup $\langle x, y\rangle$. If $x$ is of order $2^n$ for some $n > 3$, then the element $y$ can be chosen such that either $y^4=1$ and $\langle x, y\rangle$ is a generalized quaternion group with $x^y=x^{-1}$ or $y^2=1$ and $\langle x, y\rangle$ is one of the following groups: dihedral with $x^y=x^{-1}$, semidihedral with $x^y=x^{-1+2^{n-2}}$ or a group with $x^y=x^{1+2^{n-2}}$ (see \cite{Gor}, Theorem 5.4.3).  It is easy to see that from this list only generalized quaternion and dihedral subgroups can form an infinite ascending series of subgroups, so that the $2$-group $G$ can be either locally quaternion or locally dihedral, as claimed.\qed

\begin{lemma}\label{cen} Let $G$ be a group and $M$ an abelian minimal normal $p$-subgroup of $G$ for some prime $p$. Then the factor group $G/C_G(M)$ has no non-trivial finite normal $p$-subgroup.\end{lemma}

\proof Indeed, if $N/C_G(M)$ is a finite normal $p$-subgroup of $G/C_G(M)$ and $x$ is an element of order $p$ in $M$, then the subgroup $\langle N,x\rangle$ is finite and so the centralizer $C_M(N)$ is a non-trivial normal subgroup of $G$ properly contained in $M$, contrary to the minimality of $M$.\qed\medskip 

We will say that a subset $S$ of $G$ is normal in $G$ if $S^g=S$ for each $g\in G$ which means that $s^g\in S$ for every $s\in S$.

\begin{lemma}\label{DAn} Let $G$ be a group and let $A, B$ be subgroups of $G$. If a normal subset $S$ of $G$ is contained in the set $AB$ and $S^{-1}=S$, then the normal subgroup of $G$ generated by $S$ is also contained in $AB$. In particular, if $i$ is an involution with $i^G\subseteq AB$ and $N$ is the normal closure of the subgroup $\langle i\rangle$ in $G$, then $AN\cap BN=A_1B_1$ with $A_1=A\cap BN$ and $B_1=AN\cap B$. \end{lemma}

\proof Indeed, if $s,t\in S$, then $t=ab$ and $(s^{-1})^a=cd$ for some elements $a,c\in A$ and $b,d\in B$. Therefore $s^{-1}t=s^{-1}ab=a(s^{-1})^ab=(ac)(db)\in AB$ and hence the subgroup $\langle s \mid s\in S\rangle$ is contained in $AB$ and normal in $G$. Moreover, if $N$ is a normal subgroup of $G$ and $N\subseteq AB$, then it is easy to see that $AN\cap BN=(AN\cap B)N=(AN\cap B)N=(A\cap BN)(AN\cap B)$ (for details see \cite{AFG}, Lemma 1.1.3)\qed\medskip

The following slight generalization of  It\^ {o}'s theorem was proved in \cite{Sy1} (see also \cite{S_09}, Lemma 9). 

\begin{lemma}\label{genIto} Let $G$ be a group and let $A, B$ be abelian subgroups of $G$. If $H$ is a subgroup of $G$ contained in the set $AB$, then $H$ is metabelian.\end{lemma}

\section{The product of an abelian group and a group containing a quasicyclic subgroup of index $2$} 

In this section we consider groups of the form $G=AB$ with an abelian subgroup  $A$ and a subgroup  $B=X<\! y\!>$ in which $X$ is a quasicyclic subgroup of index $2$ and $y\in B\setminus X$.

\begin{lemma}\label{DAp} Let the group $G=AB$ be the product of an abelian subgroup $A$ and a non-abelian subgroup $B$ with a quasicyclic subgroup $X$ of index $2$. If $G$ has non-trivial abelian normal subgroups, then one of these is contained in the set $AX$. \end{lemma} 

\proof Suppose the contrary and let $\cal N$ be the set of all non-trivial normal subgroups of $G$ contained in the derived subgroup $G'$. Then $A_G=1$ and $ANX\ne AX$ for each $N\in\cal N$. Since $G=AB=AX\cup AXy$ and $AX\cap AXy=\emptyset$, for every $N\in\cal N$ the intersection $NX\cap AXy$ is non-empty and so $G=ANX$. Moreover, as $X=B'\leq G'$ by Lemma \ref{qua}, it follows that $G'=DNX$ with $D=A\cap G'$. It is also clear that $G=\langle A,X\rangle$, because otherwise $\langle A,X\rangle=AX$ is a normal subgroup of index $2$ in $G$. In particular, $A\cap X=1$.  

For each $N\in\cal N$ we put $A_N=A\cap BN$ and $B_N=AN\cap B$. Then $A_NN=B_NN=A_NB_N$ by \cite{AFG}, Lemma 1.1.4, and the subgroup $B_N$ is not contained in $X$, because otherwise $N$ is contained in the set $AX$, contrary to the assumption. Let $X_N=B_N\cap X$ and $C_N=A_N\cap NX_N$.  Then $X_N$ is a subgroup of index $2$ in $B_N$ and $C_NN=NX_N$ is a normal subgroup of $G$, because $(NX_N)^X=NX_N$ and $(NC_N)^A=NC_N$. Put $M=\bigcap_{N\in\cal N}N$. 

Since $G=ANX$ for each $N\in\cal N$, the factor group $G/N$ is metabelian by Lemma \ref{genIto}. Therefore also the factor group $G/M$ is metabelian and so its derived subgroup $G'/M$ is abelian. Clearly if $M=1$, then $D=A\cap G'\leq A_G=1$ and so $G'=\bigcap_{N\in\cal N}NX=X$, contrary to the assumption. Thus $M$ is an abelian minimal normal subgroup of $G$. If $M$ is finite, then its centralizer $C_G(M)$ in $G$ contains $X$ and so the group $G=A(MX)$ is metabelian. Then $G'=DMX$ is abelian and hence $D=A\cap G'\leq A_G=1$. Therefore $G'=MX$ and so $X$ is a normal subgroup of $G$,  contrary to the assumption. Thus $M$ is infinite and then $MX_M=C_MM$ is an abelian normal subgroup of $G$ with finite cyclic subgroup $X_M$ whose order is a prime power $p^k\geq 1$. 

If $M$ contains no elements of order $p$, then $M\cap C_M=1$ and $C_M^{p^k}=1$, so that $C_M$ is of order $p^k$. But then the subgroups $A_M$ and $B_M$ are of order $2p^k$ and hence the subgroup $A_MM=B_MM=A_MB_M$ is finite, contrary to the choice of $M$. Thus $M$ is an elementary abelian $p$-subgroup and so the factor group $\bar G=G/C_G(M)$ has no non-trivial finite normal $p$-subgroups by Lemma \ref{cen}. In particular, the center of $\bar G$ has no elements of order $p$. On the other hand, $G=AMX$ and $G'=DMX$ with $D=A\cap G'$, so that $\bar G=\bar A\bar X$ is a metabelian group and its derived subgroup $\bar G'=\bar D\bar X$ is abelian. Thus $\bar D\leq\bar A\cap\bar G'$ is a central subgroup of $\bar G$ which contains no elements of order $p$. But then $\bar X$ and so each of its subgroup are normal in $\bar G$. This final contradiction completes the proof.\qed\medskip 

It should be noted that if in Lemma \ref{DAp} the subgroup $B$ is locally dihedral, then the group $G=AB$ is soluble by \cite{AS_13}, Theorem 1.1. Therefore the following assertion is an easy consequence of this lemma. 

\begin{corollary}\label{Dcor} If the group $G=AB$ is the product of an abelian subgroup $A$ and a locally dihedral subgroup $B$ containing a quasicyclic subgroup $X$ of index $2$, then $AX=XA$ is a metabelian subgroup of index $2$ in $G$. \end{corollary} 

\proof Indeed, let $H$ be a maximal normal subgroup of $G$ with respect to the condition $H\subseteq AX$. If $X\leq H$, then $AH=AX$ is a metabelian subgroup of index $2$ in $G$ by It\^o's theorem. In the other case the intersection $H\cap X$ is finite and hence $HX/H$ is the quasicyclic subgroup of index $2$ in $BH/H$. Since $G/H=(AH/H)(BH/H)$ is the product of the abelian subgroup $AH/H$ and the locally dihedral subgroup $BH/H$, the set $(AH/H)(HX/H)$ contains a non-trivial normal subgroup $F/H$ of $G/H$ by Lemma \ref{DAp}.  But then $F$ is a normal subgroup of $G$ which is contained in the set $AX$ and properly contains $H$. This contradiction completes the proof.\qed\medskip

In the following lemma $G=AB$ is a group with an abelian subgroup  $A$ and a locally quaternion subgroup  $B=X<\! y\!>$ in which $X$ is the quasicyclic $2$-subgroup of index $2$ and  $y$ is an element of order $4$, so that $x^y=x^{-1}$ for each $x\in X$ and $z=y^2$ is the unique involution of $B$. It turns out that in this case the conjugacy class $z^G$ of $z$ in $G$ is contained in the set $AX$.

\begin{lemma}\label{DAc} If $G=AB$ and $A\cap B=1$, then the intersection $z^A\cap AXy$ is empty.\end{lemma}

\proof Suppose the contrary and let $z^a=bxy$ for some elements $a,b\in A$ and $x\in X$. Then $b^{-1}z=(xy)^{a^{-1}}$ and from the equality $(xy)^2=z$  it follows that $(b^{-1}z)^4=1$ and $b^{-1}zb^{-1}z=z^{a^{-1}}$. Therefore $b^{-1}z^ab^{-1}=zz^a$ and hence $bz^ab=z^az$. As  $z^a=bxy$,  we have $b(bxy)b=(bxy)z$ and so $bxyb=xyz$. Thus  $(xy)^{-1}b(xy)=zb^{-1}$. Furthermore, $bxyb^{-1}=(zb^{-1})^a$, so that $bzb^{-1}=((zb^{-1})^2)^a=(xy)^{-a}b^2(xy)^a$, i.e. the elements $z$ and $b^2$ are conjugate in $G$ by the element $g=b^{-1}(xy)^{-a}$. Since $g=cd$ for some $c\in B$ and $d\in A$, we have $b^2=z^g=z^d$ and so $z=(b^2)^{d^{-1}}=b^2$, contrary to the hypothesis of the lemma. Thus  $z^A\cap AXy=\emptyset$, as desired.\qed

\begin{theorem}\label{DAsol} Let the group  $G=AB$ be the product of an abelian subgroup $A$ and a locally quaternion subgroup $B$. If $X$ is the quasicyclic subgroup of $B$, then $AX=XA$ is a metabelian subgroup of index $2$ in $G$. In particular, $G$ is soluble of derived length at most $3$.\end{theorem} 

\proof Let $Z$ be the center of $B$, $N$ the normal closure of $Z$ in $G$ and $X=B'$, so that $X$ is the quasicyclic subgroup of index $2$ in $B$. If $A\cap B\ne 1$, then $Z$ is contained in $A\cap B$ by statement 4) of Lemma \ref{qua} and so $N=Z$. Otherwise it follows from  Lemma \ref{DAn} that $N=Z^G=Z^A$ is contained in the set $AX$. Then $N$ is a metabelian normal subgroup of $G$ by Lemma \ref{genIto} and the factor group $BN/N$ is locally dihedral by statement 7) of Lemma \ref{qua}. Since the factor group $G/N=(AN/N)(BN/N)$  is the product of an abelian subgroup $AN/N$ and the locally dihedral subgroup $BN/N$, it is soluble by \cite{AS_13}, Theorem 1.1, and so the group $G$ is soluble.

Now if $X\leq N$, then $AN=AX$ is a metabelian subgroup of index $2$ in $G$ and so the derived length of $G$ does not exceed $3$. In the other case the intersection $N\cap X$ is finite and hence $NX/N$ is the quasicyclic subgroup of index $2$ in $BN/N$. Therefore $AX=XA$ by Corollary \ref{Dcor} and this completes the proof.\qed

\section{The product of a locally quaternion and a \\ generalized dihedral subgroup}

In this section we consider groups of the form $G=AB$ with a locally quaternion subgroup $A$ and a generalized dihedral subgroup $B$.  The main part  is devoted to the proof that every group $G$ of this form has a non-trivial abelian normal subgroup. In what follows $G=AB$ is a group in which $A=Q\langle c\rangle$ with a quasicyclic $2$-subgroup $Q$ of index $2$ and an element $c$ of order $4$ such that $a^c=a^{-1}$ for each $a\in Q$ and $B=X\rtimes <\! y\!>$ with an abelian subgroup $X$ and an involution $y$  such that $x^y=x^{-1}$ for each $x\in X$.

Let $d=c^2$ denote the involution of $A$.  The following assertion is concerned with the structure of the centralizer $C_G(d)$ of $d$ in $G$. It follows from statement 4) of Lemma \ref{qua} that the normalizer of every non-trivial  normal subgroup of $A$ is contained in  $C_G(d)$.

\begin{lemma}\label{cent} The centralizer $C_G(d)$ is soluble.\end{lemma}

\proof Indeed,  if $Z=\langle d\rangle$, then the factor group $C_G(d)/Z=(A/Z)(C_B(d)Z/Z)$ is a product of the generalized dihedral  subgroup $A/Z$ and the subgroup $C_B(d)Z/Z$ which is either abelian or generalized dihedral. Therefore $C_G(d)/Z$ and thus $C_G(d)$ is a soluble group by \cite{AS_13}, Theorem 1.1, as claimed.\qed\medskip

The following lemma shows that if $G$ has no non-trivial abelian normal subgroup, then the index of $A$ in $C_G(d)$ does not exceed $2$.

\begin{lemma}\label{DAb} If $C_B(d)\ne 1$, then either $C_X(d)=1$ or $G$ contains a non-trivial abelian normal subgroup. \end{lemma}

\proof If $X_1=C_X(d)$, then $X_1$ is a normal subgroup of $B$ and $C_G(d)=AC_B(d)$. Therefore the normal closure $N=X_1^G$ is contained in $C_G(d)$,  because  $X_1^G=X_1^{BA}=X_1^A$. Since $C_G(d)$ and so $N$ is a soluble subgroup by Lemma \ref{cent}, this completes the proof.\qed\medskip 

Consider now the normalizers in $A$ of non-trivial normal subgroups of $B$. 

\begin{lemma}\label{DA} Let $G$ have no non-trivial abelian normal subgroup. If $U$ is a non-trivial normal subgroup of $B$, then $N_A(U)=1$. In particular, $A\cap B=1$. \end{lemma}

\proof Indeed, if $N_A(U)\ne 1$, then $d\in N_A(U)$ and so the normal closure ${\langle d\rangle}^G={\langle d\rangle}^B$ is contained in the normalizer $N_G(U)=N_A(U)B$. Since $N_A(U)\ne A$, the subgroup $N_A(U)$ is either finite or quasicyclic, so that $N_G(U)$ and thus ${\langle d\rangle}^G$ is soluble. This contradiction completes the proof. \qed\medskip

\begin{lemma}\label{DAk} If $C_X(d)=1$, then $G$ contains a non-trivial abelian normal subgroup.\end{lemma}

\proof Since $G=AB$, for each $x\in B$ there exist elements $a\in A$ and $b\in B$ such that $d^x=ab$. If $b\notin X$, then $b=a^{-1}d^x$ is an element of order $2$ and so $d^xad^x=a^{-1}$. As  $a^{2^k}=d$ for some $k\geq 0$, it follows that $d^xdd^x=d$ and hence $ab=d^x=(d^x)^d=(ab)^d=ab^d$. Therefore $b^d=b$ and so $b\in C_B(d)$. In particular, if $C_B(d)=1$, then $b\in X$, so that in this case the conjugacy class $d^G=d^B$ is contained in the set $AX$. 

Assume that $C_B(d)\ne 1$ and the group $G$ has no non-trivial normal subgroup. Then $C_X(d)=1$ by Lemma \ref{DAb} and without loss of generality $C_B(d)=\langle y\rangle$. Then $G=(A\langle y\rangle)X$ and so the quasicyclic subgroup $Q$ of $A$ is normalized by $y$. In particular, $d^y=d$ and the subgroup $Q\langle y\rangle$ can be either abelian or locally dihedral. We consider first the case when $y$ centralizes $Q$ and show that in this case the conjugacy class $d^G$ is also contained in the set $AX$. 

Indeed, otherwise there exist elements $a\in A$ and $b, x\in B$ such that $d^x=ab$ and $b\notin X$. Then $b\in C_B(d)=\langle y\rangle$ by what was proved above, so that $b=y$ and $d^x=ay$. As $d^B=d^{{\langle y\rangle}X}=d^X$, we may suppose that $x\in X$. But then $d^{x^{-1}}=(d^x)^y=ay=d^x$ and hence $d^{x^2}=d$. Therefore $x^2\in\langle y\rangle$ and so $x^2=1$. In particular, if $X$ has no involution, then $d^G=d^X\subseteq AX$. We show next that the case with an involution $x\in X$ cannot appear.

Clearly in this case $x$ is a central involution in $B$ and so the subgroup $D=\langle d,x\rangle$ generated by the involutions $d$ and $x$ is dihedral. It is easy to see that $d$ and $x$ cannot be conjugate in $G$ and the center of $D$ is trivial, because otherwise the centralizer $C_G(x)$ properly contains $B$, contradicting Lemma \ref{DA}. Thus $dx$ is an element of infinite order and so $D=\langle dx\rangle\rtimes\langle x\rangle$ has no automorphism of finite order more than $2$. On the other hand, if $u\in A$ and $v\in B$, then  $D^{uv}=D$ if and only if $D^u=D$ and $D^v=D$, so that $N_G(D)=N_A(D)N_B(D)$. Therefore $N_A(D)=\langle d\rangle$ and hence $z=(dx)^2$ is an element of infinite order in $N_B(D)$. But then $z\in X$ and so $\langle z\rangle$ is a normal subgroup of $B$ normalized by $d$, again contradicting Lemma \ref{DA}. Thus $X$ has no involution, as claimed.

Finally, if $N$ is the normal closure of the subgroup $\langle d\rangle$ in $G$, then $AN=NX=A_1X_1$ with $A_1=A\cap NX$ and $X_1=AN\cap X$ by Lemma \ref{DAn}. Therefore the subgroup $A_1X_1$ is soluble by Theorem \ref{DAsol}, so that $N$ and hence $G$ has a non-trivial abelian normal subgroup, contrary to our assumption. 

Thus the subgroup $Q\langle y\rangle$ is locally dihedral and so $y$ inverts the elements of $Q$. Since $A=Q\langle c\rangle$ with $a^c=a^{-1}$ for all $a\in A$, the element $cy$ centralizes $Q$ and hence the subgroup $Q\langle cy\rangle$ is abelian. But then the group $G=(Q\langle cy\rangle)B$ as the product of an abelian and a generalized dihedral subgroup is soluble by \cite{AS_13}, Theorem 1.1. This final contradiction completes the proof.\qed 

\begin{theorem}\label{QDsol} Let the group  $G=AB$ be the product of a locally quaternion subgroup $A$ and a generalized dihedral subgroup $B$. Then $G$ is soluble. Moreover, if $B$ has a quasicyclic subgroup of index $2$, then $G$ is metabelian.\end{theorem} 

\proof If $A\cap X\ne 1$, then the centralizer $C_G(d)$ is of index at most $2$ in $G$ and so $G$ is soluble by Lemma \ref{cent}.  Let $N$ be a normal subgroup of $G$ maximal with respect to the condition $A\cap NX=1$. Then $BN=(A\cap BN)B$ and the subgroup $A\cap BN$ is of order at most $2$. Therefore the subgroup $N$ is soluble and the factor group  $G/N=(AN/N)(BN/N)$ is the product of the locally quaternion subgroup $AN/N$ and the subgroup $BN/N$ which is either abelian or generalized dihedral. Hence it follows from Theorem \ref{DAsol} and Lemmas \ref{DAb} and \ref{DAk} that $G/N$ has a non-trivial abelian normal subgroup $M/N$. Put $L=MQ\cap MX$, $Q_1=Q\cap MX$ and $X_1=MQ\cap X$. Then $L=MQ_1=MX_1$ and $Q_1\ne 1$, because $A\cap MX\ne 1$ by the choice of $M$. It is also clear that $L$  is a soluble normal subgroup of $G$, because $(MQ_1)^A=MQ_1$ and $(MX_1)^B=MX_1$. Therefore the factor group $G/L$ and so the group $G$ is soluble if $AL/L$ is of order $2$. In the other case $AL/L$ is locally dihedral and $BL/L$ is abelian or generalized dihedral. Since $G/L=(AL/L)(BL/L)$, it follows that $G/L$ and so $G$ is soluble by \cite{AS_13}, Theorem 1.1. Moreover, if the subgroup $X$ is quasicyclic, then the subgroups $Q$ and $X$ centralize each other by \cite{AFG}, Corollary 3.2.8, so that $QX$ is an abelian normal subgroup of index $4$ in $G$ and thus $G$ is metabelian. \qed\medskip

The {\em Proof of Theorem} \ref{Main} is completed by a direct application of Corollary \ref{Dcor}, Theorem \ref{DAsol} and  Theorem \ref{QDsol}.\qed

\vskip 1cm 

\noindent Address of the authors:\\

\vbox{\halign to \hsize{\hbox to 0.45\hsize{#\hfil} & \hbox to
0.5\hsize{#\hfil}\cr

Bernhard Amberg           & Yaroslav P. Sysak \cr

Institut f\"ur Mathematik & Institute of Mathematics\cr

der Universit\"at Mainz & Ukrainian National Academy of Sciences\cr

D-55099 Mainz              & 01601 Kiev\cr

Germany                        & Ukraine \cr}

}

\end{document}